\documentclass[12pt,leqno]{amsart} 
\usepackage[backgroundcolor=yellow,colorinlistoftodos]{todonotes}
\usepackage{fourier}
\usepackage{datetime}
\usepackage{graphicx,float,tabularx}
\usepackage[utf8]{inputenc}
\usepackage{enumitem}   
\usepackage{amsmath,amssymb,amsfonts,mathrsfs,mathtools}
\usepackage{hyperref}
\usepackage[all,2cell]{xy}
\xyoption{arrow}
\SelectTips{cm}{10}

\newdateformat{myformat}{\THEDAY\ \monthname[\THEMONTH] \THEYEAR,
  \currenttime}
\hoffset - 2cm
\makeatletter
\renewcommand\@makefnmark{%
  \hbox{\textsuperscript{\normalfont\color{red}\@thefnmark}}}
\makeatother

\makeatletter
\providecommand\@dotsep{5}
\def\listtodoname{List of Todos}
\def\listoftodos{\@starttoc{tdo}\listtodoname}
\makeatother

\newcommand{\bookmarktodo}[2][]{\global\advance\count66 by 1
\footnote{\color{red}\textbf{#2}}%
  \addcontentsline{tdo}{todo}{\protect{
    {\color{red}[\the\count66]}\hskip 1em
    \parbox{14cm}{#2}%
  }}
}

\catcode`\@=11 \def\today{November 6, 2026}
\def\@evenfoot{\rule{0pt}{20pt}[\myformat\today] \hfill [{\tt \jobname.tex}]}
\def\@oddfoot{\rule{0pt}{20pt}{[\tt \jobname.tex}]\hfill [\myformat\today]}
\catcode`\@=13

\usepackage[object=vectorian]{pgfornament}

\usepackage{bbold,euscript}

\usepackage{quiver}

\usepackage{tikz}
\usetikzlibrary{arrows.meta}
\usetikzlibrary{decorations.markings}
\usepackage{tikz-cd}
\newtheorem{theorem}{Theorem}

\newtheorem{lemma}[theorem]{Lemma}

\theoremstyle{definition}

\newtheorem{definition}[theorem]{Definition}

\def\dn{d^{\mathbb N}}

\def\udec{{\rm dec}^{\top}}
\def\ldec{{\rm dec}_{\bot}}
\def\uldec{{\rm dec}_{\bot}^{\top}}

\def\redukce#1{\vbox to .3em{\vss\hbox{#1}}}

\def\ttN{{\mathbb N}}

\def\ttn{{\tt n}}

\def\id{{\mathbb 1}}

\def\dtop{d_{\scriptscriptstyle\top}}
\def\atop{\alpha_{\scriptscriptstyle\top}}
\def\abot{\alpha_{\scriptscriptstyle\bot}}
\def\amix{\alpha_{\scriptscriptstyle{\it{m}}}}
\def\dttop{d_{\scriptscriptstyle{\top\top}}}
\def\sttop{s_{\scriptscriptstyle{\top\top}}}
\def\stop{s_{\scriptscriptstyle\top}}
\def\dttop{d_{\scriptscriptstyle{\top\top}}}
\def\dbot{d_{\scriptscriptstyle\bot}}
\def\dbbot{d_{\scriptscriptstyle\bot\bot}}
\def\sbot{s_{\scriptscriptstyle\bot}}
\def\sbbot{s_{\scriptscriptstyle\bot\bot}}
\def\ttt#1#2{%
  \ifnum#1=0pointed
    \ifnum#2=0
    \else
      \ifnum#2=1
        {\mathbb T}_1%
      \else
        {\mathbb T}_1^{#2}%
      \fi
    \fi
  \else
    \ifnum#2=0
      \ifnum#1=1
        {\mathbb T}_2%
      \else
        {\mathbb T}_2^{#1}%
      \fi
    \else
      \ifnum#1=1
        {\mathbb T}_2%
      \else
        {\mathbb T}_2^{#1}%
      \fi
      \ifnum#2=1
        {\mathbb T}_1%
      \else
        {\mathbb T}_1^{#2}%
      \fi
    \fi
  \fi
}

\newcommand{\xrightarrowdbl}[2][]{%
  \xrightarrow[#1]{#2}\mathrel{\mkern-14mu}\rightarrow
}

\newcommand{\C}{{\tt C}}

\newcommand{\Cat}{{\tt Cat}}

\newcommand{\uu}{\mathbb{1}}

\newcommand{\NN}{\mathbb{N}}

\newcommand{\x}{\times}

\newcommand{\jed}{\mathbf{1}}

\DeclareMathOperator{\op}{op}

\DeclareMathOperator{\coker}{coker}
\DeclareMathOperator{\im}{im}

\renewcommand{\epsilon}{\varepsilon}
\renewcommand{\phi}{\varphi}
\renewcommand{\rho}{\varrho}

\title[Abelian categories and double d\'ecalage]
{Abelian categories and double d\'ecalage}
\author[M.\ Markl, D.\ Trnka]{Martin Markl and Dominik Trnka}

\catcode`\@=11
\email{markl@math.cas.cz; trnka@fme.vutbr.cz}
\@namedef{subjclassname@2010}{%
  \textup{2020} Mathematics Subject Classification}
\catcode`\@=13


\begin{document}

\begin{abstract}
We show that an additive category is abelian if and
only if its functorial nerve is a~double d\'ecalage of
a~pseudo-simplicial object in $\Cat$. In this characterization, the
canonical isomorphism between the image and coimage is
encoded by a coherent pseudo-simplicial identity.
\end{abstract}

\thanks{The authors supported by  RVO: 67985840.}
\address{Institute of Mathematics, The Czech Academy of Sciences, \v
  Zitna 25, 110 00 Praha, Czech Republic}

\maketitle

\baselineskip 17pt plus 1.5pt minus 1pt

\tableofcontents

\section*{Introduction}

Abelian categories~\cite{Freyd,maclane} are fundamental in modern
mathematics, offering an abstract setting for exactness and thus
enabling the systematic study of homology, derived functors, and
derived categories. They provide a unifying framework for structures arising
in diverse areas such as algebraic topology, representation theory,
and algebraic geometry, as well as  powerful tools and results for
these fields.

For this reason, abelian categories have been extensively studied,
leading to fundamental results such as the Freyd–Mitchell embedding
theorem~\cite{Freyd}, which states that every small abelian category 
admits a fully faithful exact embedding into a category of modules over a
ring. Considerable attention has also been devoted to the definition
itself, resulting in several equivalent characterizations in terms of
kernels and cokernels, normal monomorphisms and epimorphisms, and the
relationship between images and coimages. This last characterization
is our starting point:

\noindent\textit{An additive category $\C$ is abelian if every
  morphism $f\colon a\to b$ in $\C$ has a kernel and a cokernel, and
  the unique induced map
\begin{equation}
\label{abelian isomorphism}
  \mathrm{coim} \  f \cong \coker \ (\ker f \stackrel{\kappa_f}\longrightarrow    a)
\stackrel{\exists!} \longrightarrow \ker \ (b \stackrel{\pi_f}\longrightarrow 
\coker f) \cong  \im \  f
\end{equation}} 
\textit{is invertible.}

Our recent work \cite{abelian} provided a new characterization in terms of
lax and colax algebras over \hbox{$2$-monads}, and the present work continues
this line of research by introducing a characterization using
(pseudo)simplicial language in Theorem~\ref{hlavni}.

Let us briefly explain the results and the content of this paper. In
\cite{pasticio}, categories with kernels were characterized
simplicially by the equation $\NN_\bullet(\C) \cong \udec( X)_\bullet$, 
meaning that the
(categorical) nerve of $\C$ is the upper d\'ecalage of an almost strict
pseudo-simplicial object $X_\bullet$. In effect, this equation equips the nerve
with a (coherent) family of new functors
$d_{n+1}\colon \NN_n(\C)\to \NN_{n-1}(\C)$, which determine the
kernels via $\ker=d_2\colon \NN_1(\C)\to \NN_{0}(\C)$.

Dually, cokernels are characterized by the equation
$\NN(\C) \cong \ldec X$, which yields new functors
$d_{-1}\colon \NN_n(\C)\to \NN_{n-1}(\C)$. Combining the two, we
arrive at the double decalage construction $\uldec(-)$ (upper and lower),
which explains our title. As we show, the equation
$\NN(\C) \cong \uldec X$ equips the category~$\C$ with both kernels and
cokernels, and the isomorphism in~\eqref{abelian isomorphism} is
encoded precisely by the coherent (pseudo)simplicial
identity
\[
d_{-1}d_{3}\cong d_2d_{-1}\colon \NN_2(C)\longrightarrow\NN_0(\C).
\]

Our results emerged from the theory of operadic categories---they
complete the observations made in~\cite{bivar}. Furthermore, the
d\'ecalage construction and the pseudosimplicial setting align with
the recent fundamental work~\cite{OpCatsPs} on simplicial
characterization of the category of operadic
categories, inspired by~\cite{GarnerKockWeber}.

\vskip .5em
\noindent 
{\bf Plan of the paper.} 
In Section~\ref{Superkoherence je pro supermany.}, 
we recall supercoherent structures and introduce the
upper, lower, and double decalage constructions. Section~\ref{S Jarkou na vystavu do muzea literatury - to bude nuda!} describes
the functorial nerve of a category and its basic properties.
In Section~\ref{Jeste musime na chalupu strihat vrbicky.}, 
we recall the notions of
kernels and cokernels and their simplicial characterizations in terms
of upper and lower decalage, respectively. Finally,
Section~\ref{Pozitri zacina moje konfernce.} 
establishes the main theorem, characterizing abelian categories by
means of the double d\'ecalage construction.

\vskip .5em
\noindent 
{\bf Conventions.} We denote by $f \cdot g$ the morphism $g$ followed
by $f$, that is, the composite $f(g)$. We use the same convention for
the composites of functors and functors pre- or postcomposed with a
morphism.
The symbol $\id$, with or without
a subscript, denotes---depending on the context---the identity
endomorphism, functor, or transformation.

\section{Supercoherent structures, d\'ecalages, \&c.}
\label{Superkoherence je pro supermany.}

Here
we recall Jardine's supercoherent structures~\cite{Jardine}, 
their variants, 
and upper, lower and double d\'ecalage constructions.

\begin{definition}
    A \textit{pseudo-simplicial category} is a pseudofunctor $X\colon\Delta^{\op}\to \Cat$.
\end{definition}
According to \cite{Jardine}, a pseudo-simplicial category is generated
by the following pieces of data, which constitute a coherent 
higher-dimensional analogue to simplicial identities.
\begin{definition}
\label{definition:almost strict supercoherent}
A \textit{supercoherent structure} $X_\bullet$ consists of
categories $X_n$, $n \geq 0$, functors $$d_i : X_n\to X_{n-1}, \ 0 \leq
i \leq n, \ n \geq 1,$$
$$s_j : X_n \to X_{n+1}, \ 0 \leq j \leq n,$$
and natural isomorphisms
\begin{subequations}
\label{Jaruska upekla vybornou buchtu.}
\label{equation:simplicial identities}
\begin{align}
\label{hromada uhli}
d_id_j  &\Longrightarrow d_{j-1}d_i \ \ \text{if } i < j,
\\
\label{simpl ids with s} s_is_j  &\Longrightarrow s_{j+1}s_{i} \ \ \text{if } i \leq j, \ \hbox { and}
\\
\label{simpl ids combined} d_i s_j &\Longrightarrow
  \begin{cases*}
s_{j-1}  d_i   &  if $i < j$,  \\
\uu & if $i = j, i=j+1,$ \\ 
      s_j d_{i-1} &  if $i > j+1$. 
\end{cases*}
\end{align}
\end{subequations}
subject to coherence conditions (1.4.1)--(1.4.17) of \cite{Jardine}.
\end{definition}
Note that arrows go in the opposite direction as in \cite{Jardine}. In
particular we assume the orientation of the top-bottom simplicial
isomorphism as 
\begin{equation}
\label{Koupani v Dalesicke prehrade - zavidim.}
\amix\colon d_0d_n \Longrightarrow d_{n-1}d_0,
\end{equation}
which is a transformation of functors $X_n\to X_{n-2}$, $n\geq
2$. This convention is better suited for our purposes and, since all $\alpha$'s are isomorphisms, the direction does not really matter.

To improve readability we shall use $\dttop=d_n$ for the top face operators $d_n\colon X_{n}\to X_{n-1}$ and $\sttop=s_{n-1}$ for the top degeneracy operators  $s_{n-1}\colon X_{n-1}\to X_{n}$, for any $n\geq1$.
Similarly, we denote
\[\dtop=d_{n-1},\quad\dbbot=d_0, \quad \dbot=d_1, \quad \stop=s_{n-2}, \quad \sbot=s_1, \quad \sbbot=s_0. \]
The transformation $\amix$ of \eqref{Koupani v Dalesicke prehrade - zavidim.} is
thus written as 
\begin{subequations}
\begin{equation}
\label{Koupani v Dalesicke prehrade.}
\amix\colon \dbbot\dttop\Longrightarrow\dttop\dbbot.
\end{equation}
For the promised characterization we single out those supercoherent structures, for which most of the simplicial transformations remain strict.

\begin{definition}
\label{Uz se to rysuje.}
A supercoherent structure $U_\bullet$ is called \textit{upper almost
  strict} if  $U_0$ equals the terminal category $\jed$ and 
all the simplicial isomorphisms
of~\eqref{equation:simplicial identities} are identities, except the
top ones
\begin{equation}
\label{Posmourno.}
  \atop\colon\dttop\dttop\Longrightarrow\dttop\dtop\end{equation}
which are transformations between the functors $U_n \to U_{n-2}$, $n \geq 2$.

A supercoherent structure $L_\bullet$ is called \textit{lower almost
  strict} if $L_0 = \jed$ and if all the simplicial isomorphisms
of~\eqref{equation:simplicial identities} are identities, except the
bottom ones
\begin{equation}
\label{Dopijim kafe a pojedu za Martinem.}
\abot\colon \dbbot\dbot\Longrightarrow\dbbot\dbbot
\end{equation}
which are transformations between the functors $L_n \to L_{n-2}$, $n \geq 2$.

A supercoherent structure $X_\bullet$ is called \textit{almost strict}
if $X_0=X_1=\jed$, and all the simplicial isomorphisms
of~\eqref{equation:simplicial identities} are identities, except three
families of coherent simplicial isomorphisms \eqref{Koupani v
  Dalesicke prehrade.}, \eqref{Posmourno.}, and \eqref{Dopijim kafe a
  pojedu za Martinem.}.
 \end{definition}
\end{subequations}

For almost strict structures some of the coherence conditions
 (1.4.1)--(1.4.17) of \cite{Jardine} become trivial and the rest
 simplifies to the following two groups of equalities. 
The group (i)--(vii) relates
 $2$-cells between functors $X_{n+3}\to X_{n}$, for $n\geq 0$:

\noindent\begin{minipage}{.5\textwidth}
        \begin{enumerate}\label{10 coherence conditions}
    \item[(i)] $\abot\dbbot\cdot \abot d_2=\dbbot\abot\cdot \abot\dbot$, 
    \item[(ii)] $\amix\dbbot\cdot \dbbot\amix \cdot \abot\dttop=\dttop\abot\cdot \amix\dbot$, 
    \item[(iii)] $\atop\dbbot\cdot \dttop\amix \cdot \amix\dttop=\amix\dtop\cdot \dbbot\atop$, 
     \item[(iv)] $\atop d_{n-2}\cdot\atop\dttop=\atop\dtop\cdot \dttop\atop$, 
     \end{enumerate}
    \end{minipage}%
    \begin{minipage}{0.5\textwidth}
      \begin{enumerate}
    \item[(v)] $d_i\abot=\abot d_{i+2},$ for $1\leq i \leq n$, 
    \item[(vi)] $d_i\amix = \amix d_{i+1},$ for $1 \leq i \leq n$, 
    \item[(vii)] $d_i\atop=\atop d_i,$ for $0\leq i \leq n$. 
    \end{enumerate}    
    \end{minipage}
    
\noindent Since $X_0=X_1=\jed$ by assumption, for
$n=0,1$ the above equalities
hold trivially. 
The group (viii)--(x) below relates 
$2$-cells between functors $X_{n+1}\to X_{n}$, for $n\geq 1$: 
\begin{enumerate}
\item[(viii)] 
$\abot\sbot=\id,\qquad \abot\sbbot=\id,\qquad \abot s_j=s_{j-2}\abot,$ 
\ for $2\leq j\leq n+1$, 
    \item[(ix)] $\amix\sttop = \id,\qquad \amix\sbbot=\id,\qquad \amix s_j=s_{j-1}\amix,$ \ for $1\leq j \leq n$,
    \item[(x)] 
$\atop\stop=\id,\qquad \atop\sttop=\id,\qquad \atop s_j= s_{j}\atop,$
\ for $0\leq j\leq n-1$.
\end{enumerate}

\begin{definition}
\label{Dnes si jdu koupit ztraceny houbarsky nuz.}
The \textit{upper d\'ecalage}  $\udec(X)_\bullet$ of a 
supercoherent structure
$X_\bullet$ is the supercoherent structure with 
\[
\udec(X)_n := X_{n+1},\ n \geq
0,
\] 
 and with the simplicial operators $d^{\top}_i$, $s^{\top}_j$ 
given by the simplicial operators of
$d_i$, $s_j$ of $X_\bullet$~as 
\begin{align*}
&d^{\top}_i: \udec(X)_n = X_{n+1} \redukce{$\stackrel {d_i}\longrightarrow$} 
X_n = \udec(X)_{n-1}, \ n
\geq 1, \ 0 \leq i \leq n,  \ \hbox { and } 
\\
&s^{\top}_j: \udec(X)_n = X_{n+1} \redukce{$\stackrel {s_j}\longrightarrow$} 
X_{n+2}  =   \udec(X)_{n+1},  \ n
\geq 0, \ 0 \leq j \leq n;
\end{align*}
the top operators $d_{n+1} : X_{n+1} \to X_n$ and   
$s_{n+1} : X_{n+1} \to X_{n+2}$ of $X_\bullet$ are forgotten. The 
coherence isomorphisms
of $\udec(X)_\bullet$ are inherited from $X_\bullet$.

    The \textit{lower d\'ecalage}  $\ldec(X)_\bullet$ of a
supercoherent structure $X_\bullet$ is the supercoherent structure with 
\[
\ldec(X)_n := X_{n+1},\ n \geq
0,
\] 
and with the simplicial operators $d^{\bot}_i$, $s^{\bot}_j$ 
given by the simplicial operators of
$d_{i+1}$, $s_{j+1}$ of $X_\bullet$~as 
\begin{align*}
&d^{\bot}_i: \ldec(X)_n = X_{n+1} \redukce{$\stackrel {d_{i+1}}\longrightarrow$} 
X_n = \ldec(X)_{n-1}, \ n
\geq 1, \ 0 \leq i \leq n,  \ \hbox { and } 
\\
&s^{\bot}_j: \ldec(X)_n = X_{n+1} \redukce{$\stackrel {s_{j+1}}\longrightarrow$} 
X_{n+2}  =   \ldec(X)_{n+1},  \ n
\geq 0, \ 0 \leq j \leq n;
\end{align*}
the bottom operators $d_{0} : X_{n+1} \to X_n$ and   
$s_{0} : X_{n+1} \to X_{n+2}$ of $X_\bullet$ are forgotten. The 
coherence isomorphisms of $\ldec(X)_\bullet$ are inherited from $X_\bullet$.

The \textit{double decalage} $\uldec (X)_\bullet$ of a supercoherent structure $X_\bullet$ is obtained by applying both upper and lower decalage to the structure $X_\bullet$. Note that the order does not matter: 
\begin{equation*}
\label{decdec}
    \udec(\ldec (X))_\bullet=\uldec( X)_\bullet=\ldec(\udec (X))_\bullet.
\end{equation*}
\end{definition}

\section{The functorial nerve of a category}
\label{S Jarkou na vystavu do muzea literatury - to bude nuda!}

In this section we introduce the functorial nerve $\ttN_\bullet (\C)$
and investigate transformations between functors acting on its components.

\begin{definition}
    Let $\underline {\tt n}$ be the 
ordered set $\{0 < 1 < \cdots < n\}$ considered as a category, and
let $\ttN_n (\C) := [\underline
\ttn \to \C]$ be the $n$th piece of the categorical version of 
the standard simplicial 
nerve of a category $\C$. The latter is
the category of chains of composable arrows
\begin{equation}
\label{Mozna odjedu o den pozdeji.}
[f_1,\ldots,f_n] :=
\xymatrix@1{\bullet \ar[r]^{f_1} & \bullet \ar[r]^{f_2} 
& \cdots& \cdots  \bullet \ar[r]^{f_{n-1}}  
& \bullet \ar[r]^{f_n} & \bullet \ .
}
\end{equation}
The morphisms of $\ttN_n (\C)$ are
commutative ladders in $\C$ 
\begin{equation}
\label{S Jaruskou jsem byl na motorkach v Nepomuku.}
\xymatrix{\bullet \ar[d]^(.4){\varphi_0} \ar[r]^{f'_1} 
& \bullet \ar[r]^{f'_2}  \ar[d]^(.4){\varphi_1} 
& \cdots& \cdots  \bullet \ar[r]^{f'_{n-1}}  
& \bullet  \ar[d]^(.4){\varphi_{n-1}} \ar[r]^{f'_n} & \bullet \ar[d]^(.4){\varphi_n} 
\\
\bullet \ar[r]^{f''_1} & \bullet \ar[r]^{f''_2} 
& \cdots& \cdots  \bullet \ar[r]^{f''_{n-1}}  & \bullet \ar[r]^{f''_n} & \bullet \ .
}
\end{equation}

\end{definition}
Each commutative ladder  of length $n \geq 2$ is determined by 
two overlapping ladders of length $n-1$; namely one without $\phi_0$
and one without $\phi_n$. This Segal-type condition is expressed via
the isomorphism
\begin{subequations}
\begin{equation}
\label{je treba zabit segala}
 {\Phi_n} :    \NN_n(\C)\stackrel\cong\longrightarrow
\NN_{n-1}(\C)\x^{d_{0},d_{n-1}}_{\NN_{n-2}(\C)}\NN_{n-1}(\C), \text{
  for }n\geq 2, 
\end{equation}
where the right hand side is the pullback of $d_0$ and $d_{n-1}$. 
It will be helpful to spell out that~$\Phi_n$ sends a morphism 
$X$ of $\NN_n(\C)$ to the pair
\[
\Phi_n(X) := (d_{n}X,d_0 X)
\] 
of morphisms of $\NN_{n-1}(\C)$. 
Conversely, given morphisms $U, V$ of  $\NN_{n-1}(\C)$  
such that $d_0 U = d_{n-1}V$, there exists a unique morphism  $W$
of $\NN_n(\C)$ such that $\Phi(W) = (U,V)$. 

If we extend the
definition of the nerve by declaring $\NN_{-1}(\C)$ to be the
terminal category, the functor in~(\ref{je treba zabit segala}) becomes
\begin{equation}
\label{Dnes budu zase silet s tim vizem.}
\Phi_1 : \NN_1(\C) \longrightarrow  \NN_0(\C) \x  \NN_0(\C)
= \C \times \C,
\end{equation}
\end{subequations}
sending the morphism represented by the diagram
\[
\xymatrix@R=2em@C=2.5em{\bullet \ar[d]_(.5){\varphi_0} \ar[r]^{u} 
& \bullet  \ar[d]^(.5){\varphi_1} 
\\
\bullet \ar[r]^{v} & \bullet 
}
\]
to the pair $(\phi_0,\phi_1)$. Notice that while $\Phi_1$ is no longer an
isomorphism, it is still faithful.

The following statement prepares the ground for Lemma~\ref{Asi po 10
  letech mam chripku!}, which forms the main ingredient of the proof of
Theorem~\ref{hlavni}.  

\begin{lemma}
\label{hovadsti cinsti komunisti}
Let $c \in \{0,1,2\}$ and let $F = \{\xi_n\}_{n \geq 0}$ be a family
whose 
member $\xi_n$
is a natural transformation between 
functors \ $\NN_{n+2}(\C) \to  \NN_{n}(\C)$ such
that 
\begin{equation}
\label{Jarka nejdriv rypakovala ale pak se ji vylet libil.}
d_{i} \cdot \xi_{n+1} = \xi_{n} \cdot d_{i+c}, \ 0 \leq i \leq n+1,\
n \geq 0.
\end{equation}
Then 
\begin{itemize}
\item [(i)]
the family $F$ is uniquely determined by $\xi_0$. If $\xi_0$ is
invertible, then all higher $\xi_n$'s are invertible too.
\item [(ii)]
Any $\xi_0$, $\xi_1$ satisfying~\eqref{Jarka nejdriv
  rypakovala ale pak se ji vylet libil.} with $n =0$ are pieces of  a
unique family $F = \{\xi_n\}_{n \geq 0}$ satisfying~\eqref{Jarka nejdriv
  rypakovala ale pak se ji vylet libil.} for all $n \geq 1$. 
\end{itemize}
\end{lemma}

\begin{proof}
Let us start with the first part by proving that, given $n \geq 0$,
$\xi_n$ determines $\xi_{n+1}$. By~(\ref{je treba zabit
  segala}) and\ by the faithfulness of $\Phi_1$ in~(\ref{Dnes budu
  zase silet s tim vizem.}), it suffices to prove that $\xi_n$
determines $d_{n+1} \cdot\xi_{n+1}$ and $d_0\cdot \xi_{n+1}$. However, this
follows from~(\ref{Jarka nejdriv rypakovala ale pak se ji vylet
  libil.}), since
\[
d_{n+1}\cdot \xi_{n+1} = \xi_n\cdot d_{n+1+c} \ \hbox { and } \
d_0 \cdot\xi_{n+1} = \xi_n\cdot d_{c}.
\]
The invertibility follows from the obvious fact that $\xi_{n+1}$ is
invertible if and only if  $d_n \cdot\xi_{n+1}$ and $d_0\cdot
\xi_{n+1}$ are invertible. 
Let us prove the second part of the lemma.

Assume that we have already constructed $\xi_n$ and $\xi_{n+1}$ satisfying~(\ref{Jarka nejdriv rypakovala ale pak se ji vylet libil.})
with  some $n\geq 0$.
By~(\ref{je treba zabit segala}), $\xi_{n+2}$ 
will be determined by transformations $\omega,\lambda$ of functors 
$\NN_{n+4}(\C) \to   \NN_{n+1}(\C)$ that fulfill
\begin{equation}
\label{Za chvili volam Martinovi jestli to nezkoncilo ve spamu.}
d_0 \cdot \omega =
d_{n+1} \cdot \lambda.
\end{equation}
Such an $\xi_{n+2}$ will then be determined
by the equalities 
\[
\omega = d_{n+2} \cdot \xi_{n+2} \ \hbox {   and } \ 
\lambda = d_0 \cdot \xi_{n+2}.
\] 
Let us take
$\omega: = \xi_{n+1} \cdot d_{n+2+c}, \ 
\lambda: = \xi_{n+1} \cdot d_{c}$. 
Then 
\begin{align*}
d_0 \cdot\omega &= d_0 \cdot \xi_{n+1} \cdot d_{n+2+c} =
\xi_n  \cdot d_c  d_{n+2+c},\ \hbox { and}
\\
d_{n+1} \cdot \lambda &= d_{n+1} \cdot \xi_{n+1} \cdot d_c 
=
\xi_n \cdot d_{n+1+c} d_c. 
\end{align*}
Since $d_c  d_{n+2+c} =  d_{n+1+c} d_c$,~(\ref{Za chvili volam
  Martinovi jestli to nezkoncilo ve spamu.}) follows. Notice that, by
construction, 
\[
d_{j} \cdot \xi_{n+2} = \xi_{n+1} \cdot d_{j+c}
\]
for $j=0$ and $j=n+2$. To verify this equation for the remaining $0 <
j < n+2$ is an easy exercise on simplicial identities.
\end{proof}

\section{Simplicial characterization of kernels and cokernels.}
\label{Jeste musime na chalupu strihat vrbicky.}

In this section we 
recall the simplicial characterizations of categories 
with kernels and cokernels in terms of upper and lower d\'ecalage, respectively.

A pointed category $\C$ is a category with a chosen object 0 which is both initial and terminal. We denote the initial map of an object $a\in \C$ by $!^a\colon 0\to a$ and the terminal map by $!_a\colon a \to 0$. A~pointed category \textit{has kernels} if every map $f$ of the category has a chosen kernel
$$\ker f \xhookrightarrow{\kappa_f} a\stackrel{f}\longrightarrow b.$$
 We also assume that $\kappa_{\id_a}=!^a$ and $\kappa_{!_a}=\id_a$, for every object $a$. Similarly, a pointed category \textit{has cokernels} if every map $f$ has a chosen cokernel 
$$a\stackrel{f}\longrightarrow b\xrightarrowdbl{\pi_f} \coker f,$$
 assuming $\pi_{\id_a}=!_a$ and $\pi_{!^a}=\id_a$, for every object $a$ of the category.
 If a pointed category has both kernels and cokernels, we additionally assume $\pi_{!_a}=\kappa_{!^a}=\id_{0}$.

\begin{theorem}[Theorem 26 of \cite{pasticio}]
\label{Jen aby si stihla vyridit pas.}
A pointed category $\C$ has kernels if and only
if there is an isomorphism $\ttN_\bullet(\C) \cong\udec(U)_\bullet$,
where $U_\bullet$ is an upper almost strict supercoherent structure. 
The kernels are given by the values of $\dttop\colon \NN_{1}(\C)\to \NN_0(\C)$.
\end{theorem}

The following definition
uses the fact that each natural transformation  $\xi : p
\Rightarrow q$ between  
functors $p,q : {\tt A} \to {\tt B}$ induces a natural 
transformation $\xi^{\op} : q^{\op}
\Rightarrow p^{\op}$ between opposite functors between opposite
categories   
with components 
\[
\xi^{\op} _a : = (\xi_a)^{\op} :
q^{\op}(a) 
\to p^{\op}(a), \ a \in {\tt A}
\]
where $\xi_a : p(a) \to
q(a)$ are the components of $\xi$. 

\begin{definition}
The \textit{opposite} of an upper almost strict supercoherent
structure $U_\bullet$ is a lower almost strict supercoherent
structure $(U^{\op})_\bullet$ given by
\[
(U^{\op})_n:= (U_n)^{\op},
\] 
the opposite category of $U_n$, the
operators in dimension $n$ are
\begin{align*}
d_i^{\op}&:= (d_{n-i})^{\op}\colon (U_n)^{\op}\longrightarrow (U_{n-1})^{\op},
\ n \geq 1, \ 0 \leq i \leq n,
\\
\quad s_j^{\op}&:= (s_{n-j})^{\op}\colon (U_n)^{\op}\longrightarrow
                 (U_{n+1})^{\op}, \ n \geq 0, \ 0 \leq j \leq n.
\end{align*}
In particular, $\dttop^{\op}=(\dbbot)^{\op}$ and
$\dbbot^{\op}=(\dttop)^{\op}$. The coherence isomorphisms for
$(U^{\op})_\bullet$ are
\[\abot^{\op}:=(\atop)^{\op}\colon
  \dbbot^{\op}\dbot^{\op}\longrightarrow \dbbot^{\op}\dbbot^{\op}.
\]
\end{definition}

Dually, the opposite of a lower almost strict supercoherent structure
$L_\bullet$ is an upper almost strict supercoherent structure
$(L^{\op})_\bullet$. 
With the above definition available, we formulate a version
of Theorem~\ref{Jen aby si stihla vyridit pas.} for cokernels.   

\begin{theorem}
\label{dual:Jen aby si stihla vyridit pas.}
A pointed category \ $\C$ has cokernels if and only
if there exists an isomorphism 
\[
\ttN_\bullet(\C) \cong\ldec(L)_\bullet
\]
with  some lower almost strict supercoherent structure $L_\bullet$. 
The cokernels are given by the values of 
$\dbbot\colon \NN_{1}(\C)\to \NN_0(\C)$.\end{theorem}

\begin{proof}
A pointed category $\C$ has cokernels if and only if $\C^{\op}$ has
kernels. By Theorem~\ref{Jen aby si stihla vyridit pas.}, this is
equivalent to the isomorphism
$\NN_\bullet(\C^{\op})\cong \udec(U)_\bullet$ with $U_\bullet$ 
an upper almost strict coherent structure. Then
$L_\bullet:=(U^{\op})_\bullet$ is a lower almost strict coherent 
structure and, clearly,
\hbox{$(L^{\op})_\bullet=U_\bullet$}. Taking into
account the isomorphisms
$(\NN(\C)^{\op})_\bullet \cong \NN_\bullet(\C^{\op})$ and
$\udec(L^{\op})_\bullet \cong (\ldec (L)^{\op})_\bullet$
we conclude that
\[
(\NN(\C)^{\op})_\bullet\cong \NN_\bullet(\C^{\op})\cong
\udec(U)_\bullet=\udec(L^{\op})_\bullet \cong (\ldec (L)^{\op})_\bullet,
\]
which gives the isomorphism $(\NN(\C)^{\op})_\bullet\cong (\ldec
(L)^{\op})_\bullet$,  
equivalent to $\NN_\bullet(\C)\cong \ldec(L)_\bullet$. 
\end{proof}

\section{Simplicial characterization of abelian categories}
\label{Pozitri zacina moje konfernce.}

This section contains Theorem~\ref{hlavni}, the main result of the
paper. It combines the results of the previous section to obtain a 
simplicial characterization of abelian categories in terms of the double d\'ecalage construction.

\begin{theorem}
\label{hlavni}
An additive category $\C$ is abelian if and only if
$\uldec (X)_\bullet\cong \NN_\bullet(\C)$, where $X_\bullet$ is an
almost strict supercoherent structure of Definition~\ref{Uz se to
  rysuje.}.
    \end{theorem}

In the proof of Theorem~\ref{hlavni}, we need two auxiliary
lemmas. The first concerns  
the naturality of induced maps between kernels and
between cokernels. Its proof is straightforward, and we leave it  to the reader. 
\begin{lemma}
\label{spani v kancelari v Praze}
Consider two adjacent commutative squares as in the diagram
\[
\xymatrix{c'\ar[r]^{g'} \ar[d]_{h_0}  &  \ar[d]_{h_1}  a'\ar[r]^{f'} 
&  \ar[d]_{h_2}  b'
\\
c'' \ar[r]^{g''} & a''\ar[r]^{f''} & b''
}
\]
in which $f'g'=0$ and $f''g''=0$. The universal property of kernels induces
morphisms
\[
K\colon\ker f' \to \ker f'',\ G'\colon c'\to \ker f', \ \hbox { and } \
G''\colon c''\to \ker f'',
\]
and it holds that $KG'=G''h_0$.
Dually, the universal property of cokernels induces morphisms 
\[
C\colon\coker g' \to \coker g'',\ F'\colon \coker g'\to b', \
\hbox { and } \ F''\colon \coker g''\to b'',
\]
and it holds that $F''C=h_2F'$.
\end{lemma}

The second auxiliary lemma is a corollary of Lemma~\ref{hovadsti
  cinsti komunisti}. It refers to the families $\atop, \amix$ and
$\abot$ introduced in~\eqref{Koupani v Dalesicke
  prehrade.},~\eqref{Posmourno.} and~\eqref{Dopijim kafe a pojedu za
  Martinem.} in Section~\ref{Superkoherence je pro supermany.},
satisfying the equalities listed on page~\pageref{10 coherence
  conditions} before Definition~\ref{Dnes si jdu koupit ztraceny
  houbarsky nuz.}.  In this lemma, we will denote, for $n \geq 0$, 
the members of the
families $\abot$, $\amix$ and~$\atop$ which are transformations of
functors $X_{n+4} \to X_{n+2}$ more specifically by $\abot^n$,
$\amix^n$ and $\atop^n$, respectively. Notice that, since $X_0=X_1=\jed$, 
$\abot^0$, $\amix^0$ and $\atop^0$ are the first nontrivial members of
the corresponding families. We will also make the crucial assumption
that $\uldec( X)_\bullet \cong \NN_\bullet(\C)$ for some category
$\C$.

\begin{lemma}
\label{Asi po 10 letech mam chripku!}
The families $\abot, \amix$ and $\atop$ are uniquely
determined by $\abot^0, \amix^0$ and $\atop^0$. The members of the
families  $\abot, \amix$ and $\atop$ are invertible if and only if 
the transformations  $\abot^0, \amix^0$ and $\atop^0$ are invertible.

The first two members $\abot^i, \amix^i$ and $\atop^i$, $i=0,1$, 
uniquely extend to the corresponding full families satisfying (v),
(vi) and (vii), respectively.
Moreover, these unique extensions
automatically satisfy all the remaining conditions (i)--(x) if and
only if these conditions are satisfied by the first two  members of 
the families. 
\end{lemma}

\begin{proof}
Condition~(v) relates, for $n \geq 0$, $\abot^n$ with $\abot^{n+1}$
via the equality
\[
d_i \abot^{n+1}=\abot^n d_{i+2},\ 1  \leq i \leq n+2,
\]
which, using the isomorphism 
$\uldec( X)_\bullet \cong \NN_\bullet(\C)$, can be
rewritten as
\[
\dn_i \abot^{n+1}=\abot^n \dn_{i+2},\ 0  \leq i \leq n+1,
\]   
where $\dn_i = d_{i+1}: \NN_n(\C) \to  \NN_{n-1}(\C)$, $0\leq i \leq n$,  
are the simplicial operators of the nerve. Comparing this
with the condition in~\eqref{Jarka nejdriv rypakovala ale pak se ji vylet
  libil.},  with $\xi_n : =
\abot^n$, $n \geq 0$, Lemma~\ref{hovadsti cinsti komunisti} shows that 
the family $\abot$ is uniquely determined by $\abot^0$,
and that $\abot^0$ together with $\abot^1$ uniquely extend to the whole family
$\abot$ satisfying~(v) as claimed. The characterization of
invertibility also follows from  Lemma~\ref{hovadsti cinsti komunisti}. 
The arguments for $\amix$ and $\atop$ are similar.   

Let us prove the last part of the lemma. Assume for instance that~(i)
is satisfied for the first two members of the family $\abot$,
that is
\[
\tag{${\rm i}_0$}
\abot^0d_0 \cdot \abot^0 d_2 = d_0 \abot^{1} \cdot \abot^0 d_1.
\] 
We need to verify that 
\[ 
\tag{${\rm i}_n$}
\abot^nd_0 \cdot \abot^n d_2 = d_0 \abot^{n+1} \cdot \abot^n d_1
\]
holds for all $n \geq 1$. Since~(${\rm i}_n$) is an equality 
of transformations between functors $X_{n+5} \to
X_{n+2}$, that is between functors 
$\NN_{n+3}(\C) \to \NN_{n}(\C)$, it is enough,
by~\eqref{je treba zabit segala} if $n \geq 2$ resp.~\eqref{Dnes budu
  zase silet s tim vizem.} if $n=1$, to compare the two sides
of~(${\rm i}_n$) after
applying $\dn_0 = d_1$ and $\dn_n = d_{n+1}$. The application
of $d_1$ to the left hand side gives 
\begin{align*}
d_1(\abot^nd_0 \cdot \abot^n d_2) &= d_1\abot^nd_0 \cdot d_1\abot^n d_2
=
\abot^{n-1} d_3d_0 \cdot \abot^{n-1} d_3 d_2
=
\abot^{n-1} d_0d_4 \cdot \abot^{n-1} d_2d_4 
\\
&=
(\abot^{n-1} d_0 \cdot \abot^{n-1} d_2)d_4 
\end{align*}
where we used~(v) and the standard simplicial identities. 
Analogously we calculate 
\[
d_1(d_0 \abot^{n+1} \cdot \abot^n d_1) = d_1d_0 \abot^{n+1} \cdot 
d_1\abot^n d_1 = (d_0 \abot^n \cdot \abot^{n-1} d_1)d_4.
\]
Similarly, we verify
\begin{align*}
d_{n+1}(\abot^nd_0 \cdot \abot^n d_2) &= (\abot^{n-1}d_0 \cdot
                                        \abot^{n-1} d_2)d_{n+4}, \
                                        \hbox { and}
\\
d_{n+1}(d_0 \abot^{n+1} \cdot \abot^n d_1) &= 
(d_0 \abot^{n} \cdot
                                             \abot^{n-1} d_1)d_{n+4},
\end{align*}
We see that, for $n \geq 1$,~(${\rm i}_{n-1}$) implies~(${\rm
  i}_{n}$), so~(${\rm i}_{n}$) holds inductively for all $n \geq0$, as
claimed. The remaining items of the list can be discussed similarly.
\end{proof}

\begin{proof}[Proof of Theorem~\ref{hlavni}] 
Assume that $\C$ is abelian, in particular, that it has kernels and
cokernels.  Then by Theorems~\ref{Jen aby si stihla vyridit pas.}
and~\ref{dual:Jen aby si stihla vyridit pas.} there is an upper almost
strict supercoherent structure $U_\bullet$ and a lower almost strict
supercoherent structure $L_\bullet$ with $\udec (U)_\bullet 
\cong \NN_\bullet(\C)$
and $\ldec (L)_\bullet \cong \NN_\bullet(\C)$, such that 
$U_{n+1}\cong \NN_n(\C)\cong L_{n+1}, n\geq 0.$ Moreover, the
simplicial operators of $U_\bullet$ and $L_\bullet$ that come
from the nerve agree, and $U_0=L_0=\jed$ by
assumption.

We construct the required almost strict 
supercoherent structure $X_\bullet$ as follows. We 
define \hbox{$X_0=X_1:=\jed$}, which already determines the simplicial
operators 
\[
s_0: X_0 \to X_1,\ d_0,d_1 : X_1 \to X_0,\
\hbox { and } \  d_0,d_1,d_2 : X_2 \to X_1,
\]
as well as the coherence isomorphisms that  involve them. The operators  
$s_0,s_1: X_1 \to X_2$ send the unique object of $X_1 = \jed$ to the
null object $0 \in X_2 = \C$.  For
\hbox{$n\geq 2$} we put $X_n:=\NN_{n-2}(\C)$, with the top and bottom
simplicial operators $\dttop, \sttop$ and $\dbbot, \sbbot$ transferred
from $U_\bullet$ and $L_\bullet$, along with the accompanying coherent
simplicial isomorphisms $\atop$, $\abot$.  

The last missing pieces are
the mixed isomorphism
$\amix\colon \dbbot\dttop\Rightarrow \dttop\dbbot$. Their construction is
divided into three steps. First, we construct the transformations
$\amix$ between functors $X_n \to X_{n-2}$ for $n=4,5$ and then
we extend them to all $n \geq 4$ using Lemma~\ref{hovadsti cinsti komunisti}.
In the second step we show that 
$\amix$'s are invertible, and in the last step  we verify the
coherences (i)--(x) listed 
before Definition~\ref{Dnes si jdu koupit ztraceny houbarsky nuz.}.

\noindent\textit{Step 1: Construction of $\amix$.} The first nontrivial
$\amix$ is a transformation between functors \hbox{$X_4\to X_2$}. 
For $x = [f,g]\in \NN_2(\C)$, the component $\amix(x)$ is the
unique morphism in $X_2=\NN_0(C) = \C$ in the top tight corner in
the diagram
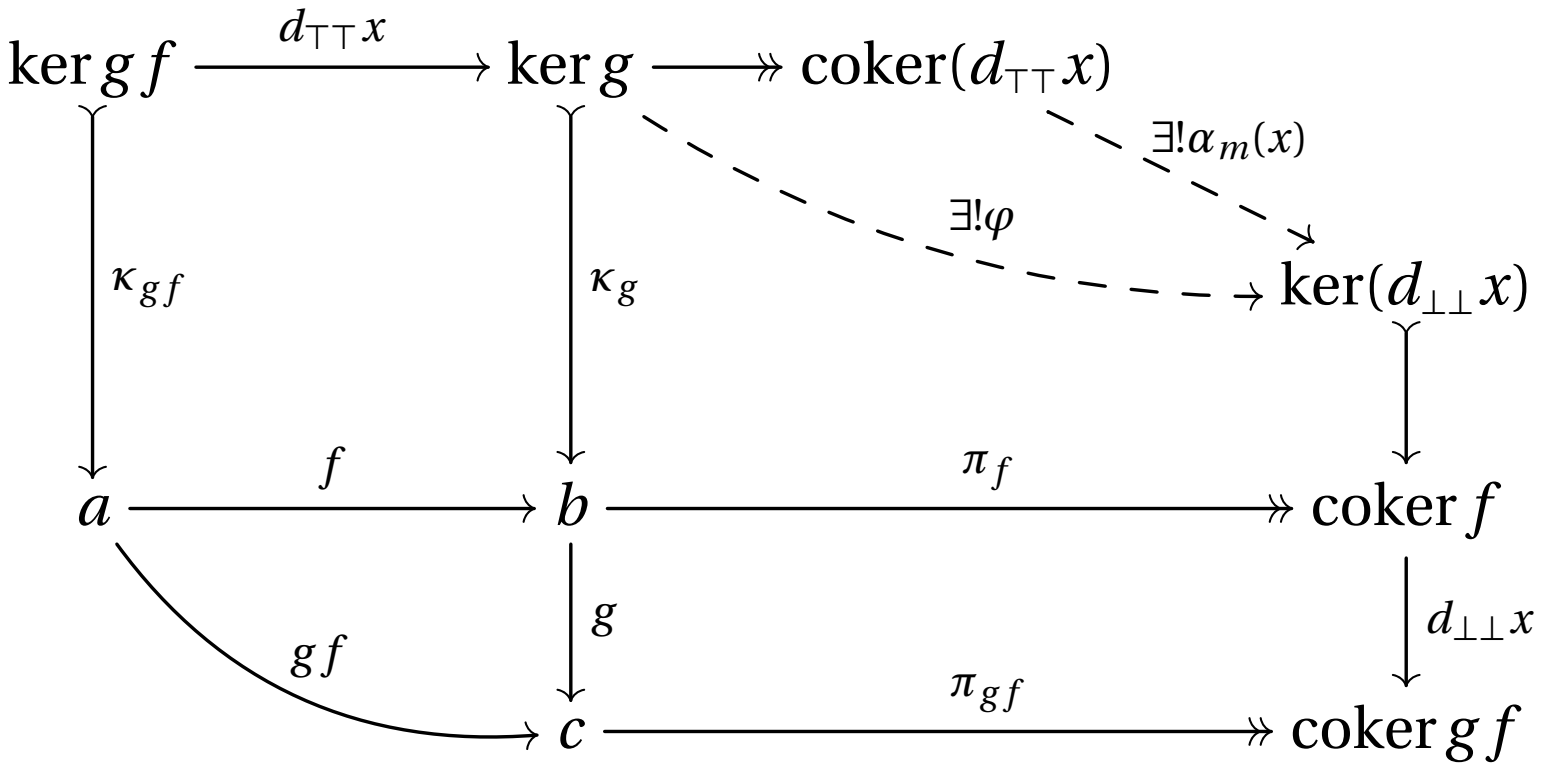
\begin{equation}
\label{tma uz o pul osme}
\xymatrix{\ker gf  \ar[rr]^{\dttop x} \ar@{>->}[dd]^{\kappa_{gf}} 
&& \ker g \ar@{->>}[r]  \ar@{>->}[dd]^{\kappa_g} \ar@{-->}@/_1.5em/[drr]^{\exists ! \phi}
& \coker (\dttop x) \ar@{-->}[dr]^{\exists ! \amix(x)}
\\
&&&& \ker  (\dbbot x) \ar@{>->}[d]
\\
a \ar@/_1.5em/[rrd]^{gf} \ar[rr]^f &&b \ar[d]^g \ar@{->>}[rr]^{\pi_f} && \coker f \ar[d]^{\dbbot x} 
\\
&& c\ar@{->>}[rr]^{\pi_{gf}} &&  \coker gf 
}
\end{equation}
in which  $\coker(\dttop
x) = \dbbot \dttop x$, $\ker (\dbbot x) = \dttop \dbbot x$ and
the dashed arrows are
induced by universal properties of kernels and cokernels. 
Let us prove the naturality of $\amix$. Namely, given a morphism  $x' \to x''$
from $x'= [f',g']$ to  $x'' =[f'',g'']$,
we must verify the commutativity of the diagram 
\begin{equation}
\label{Zitra si snad jeste naposledy v tomto roce vytahnu letadelko.}
\xymatrix@C=3.2em{
\coker (\dttop x') \ar[r]^(.53){\amix(x')} \ar[d] & 
\ker (\dbbot x') \ar[d]
\\
\coker (\dttop x'') \ar[r]^(.53){\amix(x'')} & \ker (\dbbot x'')
}
\end{equation}
in which the vertical arrows are the induced natural morphisms. To
this end, we notice the commutativity of the square
\[
\xymatrix{
\ker  g' \ar[r]^(.43){\phi'} \ar[d] & \ker (\dbbot x')\ar[d]
\\
\ker  g'' \ar[r]^(.43){\phi''} & \ker (\dbbot x'')
}
\]
which follows from the kernel version of Lemma~\ref{spani v kancelari v
  Praze} applied to
\[
\xymatrix@C=2.5em{
\ker g' \ar[r]^(.47){\pi_{f'}\kappa_{g'}} \ar[d] &  \coker f'   \ar[d]
\ar[r]^{\dbbot x'} 
&  \ar[d] \coker g'f'
\\
\ker g'' \ar[r]^(.47){\pi_{f''}\kappa_{g''}} &  \coker f'' 
\ar[r]^{\dbbot x''} 
&   \coker g''f''.
}
\]
The commutativity of~(\ref{Zitra si snad jeste naposledy v tomto roce
  vytahnu letadelko.}) 
in turn follows from the cokernel version of Lemma~\ref{spani v kancelari v
  Praze} applied to
\[
\xymatrix@C=2.5em{
\ker g'f'\ar[r]^(.54){\dttop x'} \ar[d]  &  \ar[d]   \ker g' \ar[r]^{\phi'} 
&  \ar[d] \ker (\dbbot x')
\\
\ker g''f''\ar[r]^(.54){\dttop x''}  &    \ker g'' \ar[r]^{\phi''} 
&  \ker (\dbbot x'') .
}
\]

For $n=5$, let $y\in X_5 = \NN_3(\C)$,
$y=A\xrightarrow{f}B\xrightarrow{g}C\xrightarrow{h}D$. Equality~(vi)
before Definition~\ref{Dnes si jdu koupit ztraceny houbarsky nuz.} 
gives $d_1\amix(y)=\amix(d_2y)$ and $d_2\amix(y)=\amix(d_3y)$, 
so $\amix(y)$ is represented by the diagram
\begin{equation}
    \label{diagramy nebyly zas tak hrozne}
\xymatrix@C=3em@R=2em{
Q \ar[r]^{d_0d_5(y)} \ar[d]_{\amix (d_2 y)} & R \ar[d]^{\amix (d_3 y)}
\\
O  \ar[r]^{d_4d_0(y)}  & P
}
\end{equation}
with some objects $O,P,Q,R \in \C$.
 Notice that the operators $d_2$ and $d_3$ of $X$ are the shifted
 operators $d_1^\NN$ and $d_2^\NN$ of the nerve. To show that the
 square in~\eqref{diagramy nebyly zas tak hrozne} commutes, 
we draw in Figure~\ref{vyborne tiramisu od pani z mlekarny}
the diagram describing the construction of
$\amix(y)$. 
\begin{figure}
\[
\begin{tikzpicture}[
    x=2cm,
    y=1.5cm,
    line width=0.5pt,
    line join=miter,
    line cap=butt]
    
\node[fill=white, inner sep=4pt] (A) at (2,2) {$B$};
\node[fill=white, inner sep=4pt] (B) at (3,1) {$C$};
\node[fill=white, inner sep=4pt] (C) at (3,0) {$D$};
\node[fill=white, inner sep=4pt] (D) at (0,2) {$A$};

\node[fill=white, inner sep=4pt] (E) at ($(A) + (3,0)$) {$E$};
\node[fill=white, inner sep=4pt] (F) at ($(B) + (3,0)$) {$F$};
\node[fill=white, inner sep=4pt] (G) at ($(C) + (3,0)$) {$G$};

\node[fill=white, inner sep=4pt] (H) at ($(D) + (0,2.5)$) {$H$};
\node[fill=white, inner sep=4pt] (I) at ($(A) + (0,2.5)$) {$I$};
\node[fill=white, inner sep=4pt] (J) at ($(B) + (0,2.5)$) {$J$};

\node[fill=white, inner sep=4pt] (O) at ($(E) + (0,0.7)$) {$O$};
\node[fill=white, inner sep=4pt] (P) at ($(F) + (0,0.7)$) {$P$};

\node[fill=white, inner sep=4pt] (Q) at ($(I) + (2,0)$) {$Q$};
\node[fill=white, inner sep=4pt] (R) at ($(J) + (2,0)$) {$R$};


\draw[>->,thin] (I) --  node[pos= .37,   fill=white, inner sep=2pt] {} 
node[fill=white, inner sep=2pt] {\scriptsize $\kappa_{hg}$} (A);
\draw[->>,thin] (A) --   node[pos= .31, fill=white, inner sep=2pt] {} 
 node[pos= .64, fill=white, inner sep=3pt] {}
node[fill=white, inner sep=2pt] {\scriptsize $\pi_{f}$}  (E);

\draw[->]
  (O) -- node[pos= .3, fill=white, inner sep=2pt] {\scriptsize $d_4d_0 y$} (P);
  \draw[->]
  (Q) --  node[pos= .6, fill=white, inner sep=2pt] {\scriptsize $d_0d_5 y$} (R);
\draw[thin,->]
  (R) --  node[pos= .25, fill=white, inner sep=2pt]  {\scriptsize
    $\amix d_3(y)$}  (P);defined
\draw[->]
 (Q) --  node[pos= .72, fill=white, inner sep=.2pt]  {\scriptsize
    $\amix d_2(y)$}  node[pos= .57, fill=white, inner sep=2pt] {}   (O);

\draw[thin,dashed,-<] 
  (I) ..controls ++(0.5,-2.3) and ++(-0.1,0) .. 
node[pos= .15, fill=white, inner sep=2pt] {}
node[pos= .37, fill=white, inner sep=12pt] {}
node[fill=white, inner sep=2pt] {\scriptsize $\phi'$}  (O);
\draw[dashed,-<] 
  (J) ..controls ++(0.5,-2.3) and ++(-0.1,0) .. node[fill=white, inner sep=2pt] {\scriptsize $\phi''$} (P);

 \draw[thin,<-]
  (C) .. controls ++(-0.5,0.2) and ++(0,-1) .. node[pos= .5,
  fill=white, inner sep=2pt] {\scriptsize $hg$} 
node[pos= .65, fill=white, inner sep=5pt] {}  (A);
\draw[thin,->] (A) --  node[pos= .5, fill=white, inner sep=2pt]
{\scriptsize $g$} (B);
\draw[->] (B) -- node[pos= .5, fill=white, inner sep=2pt]
{\scriptsize $h$} (C);
\draw[thin,<-] (C) .. controls ++(-1,0) and ++(1,-2) ..  node[pos= .5, fill=white, inner sep=2pt]
{\scriptsize $hgf$} (D);
\draw[thin,->] (D) -- node[pos= .5, fill=white, inner sep=2pt]
{\scriptsize $f$}  (A);

\draw[thin,<-]
  (B) .. controls ++(-1,0) and ++(1,-1) .. 
node[pos= .5, fill=white, inner sep=2pt] {\scriptsize $gf$}  (D);


\draw[thin,->] (E) -- node[pos= .3, fill=white, inner sep=2pt] {} (F);
\draw[thin,->] (F) -- (G);
\draw[thin,<-]
  (G) .. controls ++(-0.5,0.2) and ++(0,-1) .. node[pos= .65,
  fill=white, inner sep=2pt] {} 
 node[pos= .94,
  fill=white, inner sep=2pt] {}  (E);


\draw[thin,->] (H) -- node[fill=white, inner sep=2pt] {\scriptsize
  $d_5d_3 y$}  (I);
\draw[thin,->] (I) -- (J);
\draw[thin,<-]
  (J) .. controls ++(-1,0) and ++(1,-1) .. node[fill=white, inner
  sep=2pt] {\scriptsize $d_5d_2 y$}  (H);




\draw[thin,<-<] 
(D) -- node[fill=white, inner sep=2pt] {\scriptsize $\kappa_{hgf}$}  (H);

 
\draw[<<-,thin] (Q) -- (I);
\draw[<-<,thin] (E) -- (O);


\draw[<<-,thin] (R) -- (J);
\draw[>->,thin] (J) -- node[fill=white, inner sep=2pt] {\scriptsize
  $\kappa_{h}$}  (B);
\draw[->>,thin] (B) --  
node[fill=white, inner sep=2pt] {\scriptsize $\pi_{gf}$}  (F);
\draw[<-<,thin] (F) -- (P);

\draw[thin,->>] 
(C) -- node[fill=white, inner sep=2pt] {\scriptsize $\pi_{hgf}$}  (G);
\end{tikzpicture}
\]
\caption{Construction of $\amix (y)$.\label{vyborne tiramisu od pani z mlekarny}}
\end{figure}
In that diagram, $E\to F\to G=\dbbot y$, and $H\to I\to J=\dttop y$.
Applying the kernel version of Lemma~\ref{spani v kancelari v Praze}
to the diagram
\[
\xymatrix@C=3em{
I\ar[r]^{\kappa_{hg}} \ar[d] &\ar[r]^{\pi_f}  B&E \ar[r]^{d_0d_3 y} \ar[d] 
 &G \ar@{=}[d]
\\
J\ar[r]^{\kappa_h} &C \ar[r]^{\pi_{gf}}& \ar[r]^{d_0 d_2y}   F&G
}
\]
in which $I\to J$ is the induced morphism between kernels, and $E\to
F$ is the induced morphism between cokernels, we verify that the 
square $(IJOP)$ commutes. Finally, the cokernel version of
Lemma~\ref{spani v kancelari v Praze} applied to 
\[
\xymatrix@C=3em{
H\ar[r]^{d_5d_3 y}\ar@{=}[d]  & \ar[r]^{\phi'} \ar[d]  I&O \ar[d] 
\\
H  \ar[r]^{d_5d_2 y} &  \ar[r]^{\phi''} J&P
}
\]
gives the required commutativity of~\eqref{diagramy nebyly zas tak hrozne}. 
The naturality of $\amix$ on $X_5$ follows from the fact that it is
determined by $\amix(d_2)$ and $\amix(d_3)$, cf.~\eqref{diagramy
  nebyly zas tak hrozne}, and that $\amix$ is natural on $X_4$ as we
already established.

We have constructed the two initial terms of the family $\amix$.
Lemma~\ref{Asi po 10 letech mam chripku!} tells us that
they extend uniquely to the full family satisfying~(vi). 
It remains to prove that
the families $\abot$, $\amix$ and $\atop$ together satisfy all
conditions (i)--(x).
By Lemma~\ref{Asi po 10 letech mam chripku!} again, it will suffice
to check that they are satisfied by the first two members of 
the families. This will be done in Step~3 below.

\vskip .5em
\noindent
\textit{Step 2: Invertibility of $\amix$.}
By Lemma~\ref{Asi po 10 letech mam chripku!},
it suffices to show that $\amix(x)$ is an invertible morphism for 
any $x = [h,g]\in \NN_2(\C)$.
If $\C$ is an abelian category of modules over a ring $R$, this was
already observed in~\cite[Proposition~26]{bivar}. The morphisms
$h_Z$, $g^S$, and $\rho$ in diagram~(25) of \cite{bivar} are precisely
morphisms $\dttop(x)$, $\dbbot(x)$, and the component $\amix(x)$,
respectively.  

If $\C$ is arbitrary abelian, we consider its small
abelian subcategory generated by the diagram
$S\xrightarrow{h} T\xrightarrow{g}Z$. Freyd–Mitchell's embedding
theorem~\cite[Theorem~4.44 \& page~104]{Freyd}
gives an exact fully faithful embedding of this subcategory
into the category of modules, proving invertibility of $\rho$ in the
general case.

\vskip .5em
\noindent
\textit{Step 3: Coherence.} We need to verify equalities
(ii), (iii) and (ix) only. Indeed, 
equalities (i), (v) and (viii) hold because $L$ is lower 
almost strict, while  (iv), (vii) and (x) hold because $U$ is
upper almost strict. Moreover, since  
the higher $\amix$'s 
are inductively 
constructed from the lower ones, equality (vi) follows from the
second part of Lemma~\ref{Asi po 10 letech mam chripku!}. 

Condition (ii) is dual to (iii), so it suffices to prove one of
them, say (iii). 
By Lemma~\ref{Asi po 10 letech mam chripku!} again, 
it is enough to verify (iii) for the first two
members of the families, where it is an equality of  
transformations of functors $X_5 \to X_2$, that is, of functors $\NN_3(\C)
\to \NN_0(\C)$. Let therefore 
\[
y = a \xrightarrow f b \xrightarrow g c
\xrightarrow h d \in \NN_3(\C).
\]   
Equality (iii) evaluated at $y$ appears as the equality
\begin{equation}
\label{Predevcirem jsem priletel z Ciny.}
S\to T\to U  \ = \ S\to W\to V \to U
\end{equation}
of the composites of morphisms of $\C$ in the diagram in
Figure~\ref{Jarka mi koupila Tonic!}.
\begin{figure}
\begin{center}
\begin{tikzpicture}[
   x=2cm,
    y=1.5cm,
    line width=0.5pt,
    line join=miter,
    line cap=butt]
\node[fill=white, inner sep=4pt] (A) at (2,2) {$B$};
\node[fill=white, inner sep=4pt] (B) at (3,1) {$C$};
\node[fill=white, inner sep=4pt] (C) at (3,0) {$D$};
\node[fill=white, inner sep=4pt] (D) at (0,2) {$A$};

\node[fill=white, inner sep=4pt] (E) at ($(A) + (3,0)$) {$E$};
\node[fill=white, inner sep=4pt] (F) at ($(B) + (3,0)$) {$F$};
\node[fill=white, inner sep=4pt] (G) at ($(C) + (3,0)$) {$G$};

\node[fill=white, inner sep=4pt] (H) at ($(D) + (0,2.5)$) {$H$};
\node[fill=white, inner sep=4pt] (I) at ($(A) + (0,2.5)$) {$I$};
\node[fill=white, inner sep=4pt] (J) at ($(B) + (0,2.5)$) {$J$};

\node[fill=white, inner sep=4pt] (K) at ($(D) + (-1,1)$) {$K$};
\node[fill=white, inner sep=4pt] (L) at ($(H) + (-1,1)$) {$M$};
\node[fill=white, inner sep=4pt] (M) at ($(I) + (-1,1)$) {$N$};
\node[fill=white, inner sep=4pt] (N) at ($(A) + (-1,1)$) {$L$};

\node[fill=white, inner sep=4pt] (O) at ($(E) + (0,0.7)$) {$O$};
\node[fill=white, inner sep=4pt] (P) at ($(F) + (0,0.7)$) {$P$};

\node[fill=white, inner sep=4pt] (Q) at ($(I) + (2.5,0)$) {$Q$};
\node[fill=white, inner sep=4pt] (R) at ($(J) + (2.5,0)$) {$R$};

\node[fill=white, inner sep=4pt] (S) at ($(E) + (-0.5,0.5)$) {$U$};
\node[fill=white, inner sep=4pt] (T) at ($(O) + (-0.5,0.5)$) {$V$};
\node[fill=white, inner sep=4pt] (U) at ($(Q) + (-0.5,0.5)$) {$W$};

\node[fill=white, inner sep=4pt] (V) at ($(M) + (1.5,0)$) {$S$};
\node[fill=white, inner sep=4pt] (W) at ($(N) + (1.5,0)$) {$T$};


 \draw[->>] 
  (N) --
node[pos= .7, fill=white, inner sep=2pt] {}  
 (W);
\draw[>->] (I) -- 
node[pos= .2, fill=white, inner sep=2pt] {\scriptsize$\kappa_{hg}$}
node[pos= .38, fill=white, inner sep=2pt] {}
(A);   
\draw [->] (V) --              
node[pos= .2, fill=white, inner sep=2pt] {\scriptsize $\dbbot \atop y$}(W)
node[pos= .61, fill=white, inner sep=3pt] {}
node[pos= .85, fill=white, inner sep=3pt] {}
node[pos= .37, fill=white, inner sep=3pt] {};
\draw [->] (W) --               
node[pos= .5, fill=white, inner sep=2pt] {\scriptsize $\amix \dtop   y$}(S)
node[pos= .21, fill=white, inner sep=3pt] {};
\draw [->] (T) -- (S);           
\draw [->] (U) --           
node[pos= .5, fill=white, inner sep=2pt] {\scriptsize $\dttop \amix  y$}
node[pos= .21, fill=white, inner sep=3pt] {}
node[pos= .91, fill=white, inner sep=3pt] {}
(T);
\draw [->] (V) --       
node[pos= .5, fill=white, inner sep=2pt] {\scriptsize $\amix\dttop y$}
(U);


\draw[<-] (C)           
.. controls ++(-0.5,0.2) and ++(0,-1) .. 
node[pos= .45, fill=white, inner sep=2pt] {\scriptsize $hg$}
node[pos= .65, fill=white, inner sep=2pt] {}
(A);
\draw[->] (A) --   
node[pos= .5, fill=white, inner sep=2pt] {\scriptsize $g$} (B);
\draw[->] 
(B) -- 
node[pos= .5, fill=white, inner sep=2pt] {\scriptsize $h$} (C);    
\draw[<-] (C)          
.. controls ++(-1,0) and ++(1,-2) .. 
node[pos= .5, fill=white, inner sep=2pt] {\scriptsize $hgf$}
(D);
\draw[->] (D) --    
node[pos= .5, fill=white, inner sep=2pt] {\scriptsize $f$}
 (A);
\draw[<-] (B)      
.. controls ++(-1,0) and ++(1,-1) ..
node[pos= .5, fill=white, inner sep=2pt] {\scriptsize $gf$}
 (D);

\draw[->] (E) -- 
node[pos= .5, fill=white, inner sep=2pt] {\scriptsize $\dbbot\dtop y$}
(F);
\draw (F)[->] -- (G);
\draw[<-] (G) .. controls ++(-0.5,0.2) and ++(0,-1) ..
node[pos= .65, fill=white, inner sep=2pt] {}
 (E);


\draw[<-] (K) -- (L);   
\draw[->] (L) --        
node[pos= .5, fill=white, inner sep=2pt] {\scriptsize $\dttop\dttop y$}
 (M);
\draw[->] (M) -- 
node[pos= .72, fill=white, inner sep=2pt] {}
node[pos= .38, fill=white, inner sep=2pt] {}
 (N);
\draw[<-] (N) --           
node[pos= .5, fill=white, inner sep=2pt] {}
node[pos= .25, fill=white, inner sep=2pt] {\scriptsize $\dttop \dtop y$}
(K);      


\draw[->] (H) -- (I);
\draw[->] (I) -- (J);
\draw[<-] (J) .. controls ++(-1,0) and ++(1,-1) .. (H);


\draw[->] (O) --    
node[pos= .5, fill=white, inner sep=2pt] {\scriptsize $\dttop \dbbot y$}(P);
\draw[->] (R) --
 (P);
\draw[->] (Q) -- (R);
\draw[->] (Q) --
node[pos= .56, fill=white, inner sep=2pt] {} 
(O);

 \draw[>->]      
  (N) -- node[pos= .6, fill=white, inner sep=2pt]
  {\scriptsize$\kappa_{g}$} (A);
\draw[>->] (K)       
--node[pos= .4, fill=white, inner sep=2pt] {\scriptsize$\kappa_{gf}$} (D);
\draw[<-<] (D) --      
node[pos= .6, fill=white, inner sep=2pt] {\scriptsize$\kappa_{hgf}$}
 (H);
\draw[<-<] (H) -- (L);      

\draw[>->] 
(M) -- (I);      

\draw[<<-] (Q) -- (I);
\draw[->>] (A)   
 -- node[pos= .6, fill=white, inner sep=2pt] {\scriptsize$\pi_{f}$}
node[pos= .31, fill=white, inner sep=2pt] {} (E);
\draw[>->] (O) -- (E);

\draw[->>] (J) -- (R);
\draw[>->] (J)     
-- node[pos= .5, fill=white, inner sep=2pt] {\scriptsize$\kappa_{h}$} (B);
\draw[->] (B)    
--  node[pos= .5, fill=white, inner sep=2pt] {\scriptsize$\pi_{gf}$} (F);
\draw[>->] (P) -- (F);


  \draw[->>]           
  (C) -- node[pos= .5, fill=white, inner sep=2pt] {\scriptsize$\pi_{hgf}$}  (G);
  \draw[>->]   
  (U) -- (Q);
  \draw[>->]     
  (T) -- (O);
  \draw[>->] 
  (S) -- (E);   
  \draw[->>] 
  (M) -- (V);
\end{tikzpicture}
\end{center}
\caption{The coherence.\label{Jarka mi koupila Tonic!}}
\end{figure}
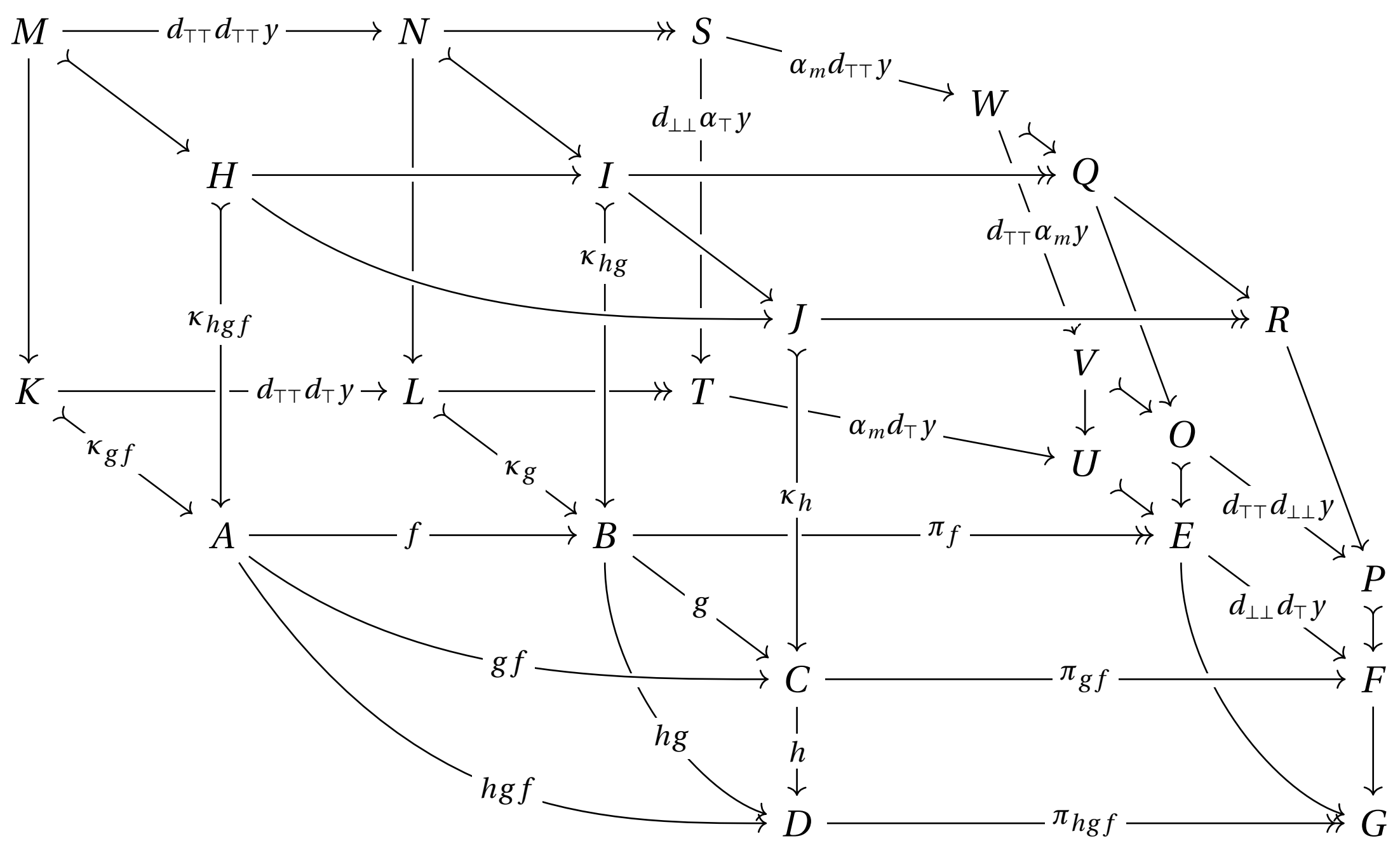
The arrow $N\to S$ is a cokernel map, hence an
epimorphism, and the arrow $U\to E$ is a kernel map,
hence a monomorphism. Therefore,~\eqref{Predevcirem jsem
  priletel z Ciny.} is
equivalent to 
\[
N\to S \to T \to U\to E  \ = \ N \to S \to W \to V \to U\to E,
\] 
which can be verified by the standard diagram
chase. Notice that all other regions of the diagram commute. Thus (iii),
and hence also (ii), have been verified.

By Lemma~\ref{Asi po 10 letech mam chripku!}, it is enough 
to prove the equalities in (ix)
only when they compare transformations of functors $X_3 \to X_2$, that
is, functors $\NN_1(\C) \to \NN_0(\C)$. Let us evaluate, for example, $
\amix\sbbot=\id$ at $f : a \to b  \in \NN_1(\C)$.
Diagram~\eqref{tma uz o pul osme} defining
$\amix x$ for $x = \sbbot f=[!^a,f]$ becomes
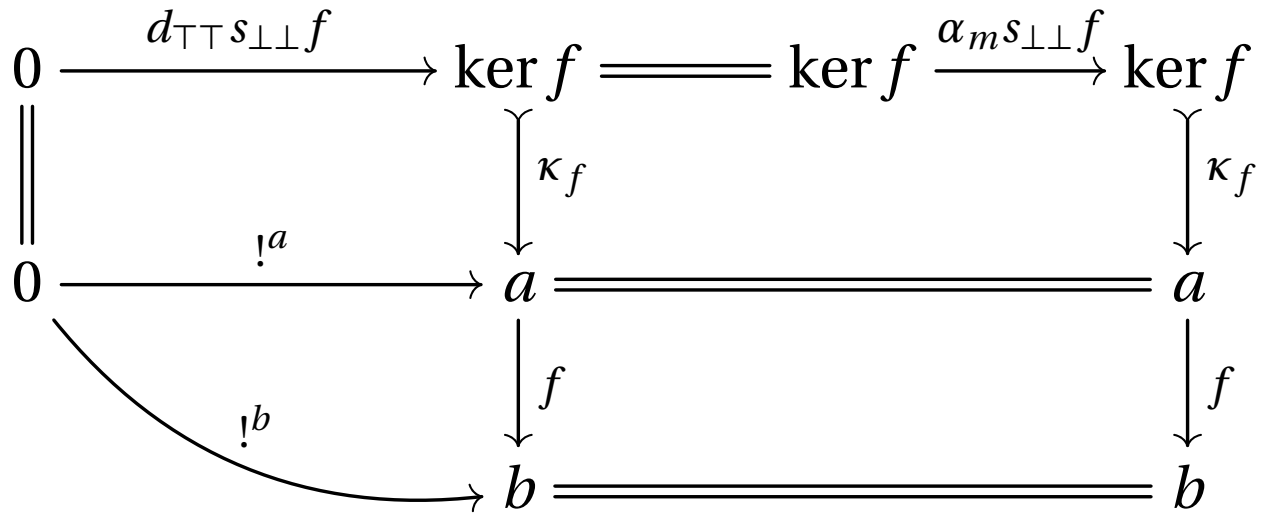
\[
\xymatrix@C=3em{
0 \ar[rr]^(.43){\dttop \sbbot f} \ar@{=}[d] 
&& \ker f \ar@{=}[r]  \ar@{>->}[d]^{\kappa_f} 
& \ker f  \ar@{->}[r]^{\amix\sbbot f}
& \ker f \ar@{>->}[d]^{\kappa_f}
\\
0 \ar@/_1.5em/[rrd]^{!^b} 
\ar[rr]^{!^a} &&a \ar[d]^f \ar@{=}[rr] && a \ar[d]^{f} 
\\
&& b\ar@{=}[rr] &&  b
}
\]
so clearly $\amix\sbbot f = \id$ as required.
The remaining equalities in (ix) can be verified analogously. This
concludes the proof that, for an abelian category $\C$,
$\NN_\bullet(\C) \cong \uldec (X)_\bullet$ for some almost strict supercoherent structure  $X_\bullet$.

\noindent\textit{Opposite implication.}
Assume that $\uldec(X)_\bullet\cong\NN_\bullet(\C)$ for an almost
strict supercoherent structure $X_\bullet$. Since then
$\NN_\bullet(\C) \cong\udec(\ldec (X))_\bullet$, Theorem~\ref{Jen aby si stihla
  vyridit pas.} implies that 
the category $\C$ has kernels. Dually, 
Theorem~\ref{dual:Jen aby si stihla vyridit pas.} implies the
existence of cokernels. It remains to establish the
isomorphism~\eqref{abelian isomorphism}.
To this end, we use the factorization  
\begin{equation}
\label{Zitra mi zacinaji prijizdet ucastnici.}
\iota_f := a \xrightarrow {\id_a \oplus f} a \oplus b
\xrightarrow{0+\id_b} b
\end{equation}
regarded as an object of $X_4 = \NN_2(\C)$, and the associated
diagram~\eqref{tma uz o pul osme} with $x = \iota_f$. As in 
the proof of~\cite[Theorem~24]{abelian}, we observe that this diagram
holds for an arbitrary choice of concrete (co)kernels in $\C$, which
leads to
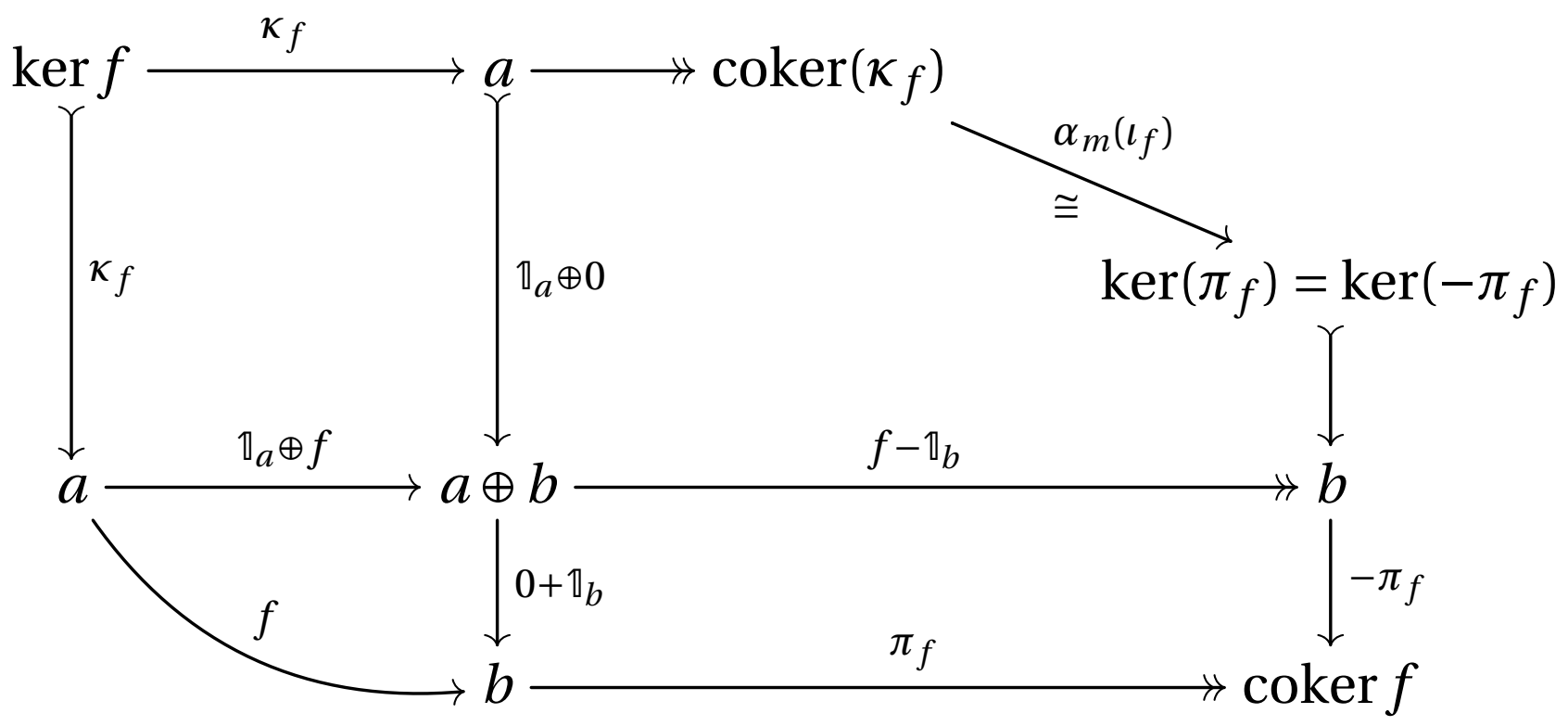
\[
\xymatrix{\ker f  \ar[rr]^{\kappa_f} \ar@{>->}[dd]^{\kappa_{f}} 
&& a\ar@{->>}[r]  \ar@{>->}[dd]^{\id_a \oplus 0}
& \coker (\kappa_f) \ar@{->}[dr]^{\amix(\iota_f)}_\cong
\\
&&&& \ker  (\pi_f) = \ker  (-\pi_f) \ar@{>->}[d]
\\
a \ar@/_1.5em/[rrd]^{f} \ar[rr]^{\id_a \oplus f} && a\oplus b 
\ar[d]^{0 + \id_b} \ar@{->>}[rr]^{f- \id_b} && b \ar[d]^{-\pi_f} 
\\
&& b\ar@{->>}[rr]^{\pi_{f}} &&  \coker f 
}
\]
in which the isomorphism $\amix(\iota_f)$ comes from the
supercoherent structure. This concludes the proof of Theorem~\ref{hlavni}.
\end{proof}

The additivity of the category $\C$ was used in the second part of the
proof of Theorem~\ref{hlavni} merely to construct the
factorization in~\eqref{Zitra mi zacinaji prijizdet ucastnici.}.
Factorizations of the same type exist under weaker assumptions,
formulated in

\begin{definition}
A category with kernels and cokernels has a
  \textit{$KC$-factorization} if every morphism $f\colon a\to b$
can be written as $f=em$ with a monomorphism $m$ and an epimorphism
$e$, such that  there are isomorphisms
$\xi\colon a\xrightarrow\cong \ker(e)$ and $\zeta\colon
\coker(m)\xrightarrow\cong b$ 
satisfying
$m \,\kappa_f=\kappa_e \,\xi \, \kappa_f$ and $\pi_f \,e=\pi_f \,\zeta\,\pi_m$.
\end{definition}

Factorization~\eqref{Zitra mi zacinaji prijizdet ucastnici.} is
clearly of the type required by the above definition.
Recall from \cite[S.~1.1.]{BorceuxGrandis} that a
Puppe-exact category is a pointed category with kernels and cokernels
such that the canonical map $\coker(\kappa_f)\to \ker(\pi_f)$ is an
isomorphism for every morphism $f$. Our theory provides the following
characterization. 
    
\begin{theorem}
    A category $\C$ with a $KC$-factorization is Puppe-exact if and
    only if there is an isomorphism $\uldec (X)_\bullet\cong
    \NN_\bullet(\C)$ for some almost strict
    supercoherent structure $X_\bullet$.
    \end{theorem}

\bibliographystyle{plain} 

\end{document}